\input amstex
\documentstyle{amsppt}
\pagewidth{5.4in}
\pageheight{7.6in}
\magnification=1200
\TagsOnRight
\NoRunningHeads
\topmatter
\title
\bf Uniqueness of solutions of Ricci flow on complete noncompact manifolds
\endtitle
\author
Shu-Yu Hsu
\endauthor
\affil
Department of Mathematics\\
National Chung Cheng University\\
168 University Road, Min-Hsiung\\
Chia-Yi 621, Taiwan, R.O.C.\\
e-mail:syhsu\@math.ccu.edu.tw
\endaffil
\date
Oct 7, 2011
\enddate
\address
e-mail address:syhsu\@math.ccu.edu.tw
\endaddress
\abstract
We give a simple proof of the uniqueness of solutions of the Ricci flow 
on complete noncompact manifolds with bounded curvatures using the De 
Turck approach. As a consequence we obtain a correct proof of the existence 
of solution of the Ricci harmonic flow on complete noncomplete manifolds 
with bounded curvatures. We give a simple example of a complete manifold
with bounded curvature and injectivity radius tending to zero as the point 
goes to infinity. We also give strong proof and argument why the crucial 
lemma Lemma 2.2 of \cite{CZ} cannot be valid.
\endabstract
\keywords
Ricci harmonic flow, Ricci flow, existence, uniqueness 
\endkeywords
\subjclass
Primary 58J35, 53C43 Secondary 35K55
\endsubjclass
\endtopmatter
\NoBlackBoxes
\define \pd#1#2{\frac{\partial #1}{\partial #2}}
\define \1{\partial}
\define \2{\overline}
\define \3{\varepsilon}
\define \4{\widetilde}
\define \5{\underline}
\define \R{\Bbb{R}}
\define \ov#1#2{\overset{#1}\to{#2}}
\define \oa#1{\overset{a}\to{#1}}
\document

Recently there is a lot of study on the Ricci flow on manifolds by 
R.~Hamilton \cite{H1--3} and others. 
Existence of solution $(M,g(t))$, $0\le t\le T$, of the Ricci flow equation 
$$ 
\frac{\1 }{\1 t}g_{ij}=-2R_{ij}\tag 0.1
$$
on compact manifold $M$ where $R_{ij}(t)$ is the Ricci curvature of $g(t)$ 
and $g_{ij}(x,0)=g_{ij}(x)$ is a smooth metric on $M$ is proved by 
R.~Hamilton in \cite{H1}. R.~Hamilton \cite{H1} also proved that when 
$g_{ij}(x)$ is a metric of strictly positive Ricci curvature,
then the evolving metric will converge modulo scaling to a metric 
of constant positive curvature. 

Since the proof of existence of solution of the Ricci flow in \cite{H1} is
very hard, later D.M.~DeTurck \cite{D} deviced another method to prove 
existence and uniqueness of solution of Ricci flow on compact manifolds. 
Let $M$ be a n-dimensional manifold with $(M,g_{ij}(t))$, $0\le t\le T$, 
being a solution of the Ricci flow (0.1) 
and let $(N,h_{\alpha\beta})$ be a fixed $n$-dimensional manifold. He 
introduced the associated Ricci harmonic flow $F=(F^{\alpha}):
(M,g(t))\to (N,h)$ given by
$$
\frac{\1F}{\1 t}=\Delta_{g(t),h}F\tag 0.2
$$
where
$$
\Delta_{g(t),h}F^{\alpha}=\Delta_{g(t)}F^{\alpha}+g^{ij}(x,t)
\4{\Gamma}^{\alpha}_{\beta ,\gamma}(F(x,t))\frac{\1F^\beta}{\1 x^i}
\frac{\1F^\gamma}{\1 x^j}\tag 0.3
$$
in the local co-ordinates $x=(x^1,\dots,x^n)$ of the domain manifold
$(M,g_{ij}(t))$ and the local co-ordinates $(y^{\alpha})$ of the target 
manifold $(N,h_{\alpha\beta})$ with 
$$
\Delta_{g(t)}F^{\alpha}=g^{ij}\nabla_i\nabla_jF^{\alpha}
$$
and $\4{\Gamma}^{\alpha}_{\beta ,\gamma}$ 
being the Christoffel symbols of $(N,h_{\alpha\beta})$. 
When the solution $F(\cdot ,t)$ of (0.2) is a diffeomorphism, 
this harmonic map flow induces a push forward metric 
$$
\hat{g}(t)=(F)_{\ast}(g(t))=(F(\cdot,t)^{-1})^{\ast}(g(t))\tag 0.4
$$
on the target manifold $N$ which satisfies the Ricci-DeTurck flow \cite{H3},
$$
\frac{\1 }{\1 t}\hat{g}_{\alpha\beta}=(L_V\hat{g})_{\alpha\beta}
-2\hat{R}_{\alpha\beta}\tag 0.5
$$
for some time varying vector field $V$ on the target manifold $N$ 
where $\hat{R}_{\alpha\beta}$ is the Ricci curvature associated with the 
metric $\hat{g}(t)$. Since (0.5) is strictly parabolic \cite{H3}, it is
easiler to solve (0.5) than (0.1) which is weakly parabolic \cite{H1}. 
The existence 
and uniqueness of solutions of Ricci flow on compact manifolds are then 
reduced to the study of existence and other properties of the harmonic map 
flow (0.2) and the Ricci-DeTurck flow (0.5). We refer the reader to the paper 
\cite{H3} of R.~Hamilton on a sketch of this approach on compact manifolds. 

Naturally one would expect this approach should also work for non-compact
complete Riemmanian manifolds. In \cite{S1} W.X.~Shi used this technique
to prove the existence of solution of (0.1) on complete non-compact Riemannian
manifolds. In [LT] P. Lu and G.~Tian used the De Turck trick to prove the
uniqueness of the standard solution of Ricci flow on $\Bbb{R}^n$, $n\ge 3$,
which is radially symmetric about the origin. Recently S.Y.~Hsu \cite{Hs}
extended the result of \cite{LT} and proved the uniqueness of the solution 
of the radially symmetric solution of the Ricci harmonic flow (0.2) 
assoicated with the standard solution of Ricci flow. 

In \cite{CZ} B.L.~Chen and X.P.~Zhu attempted to prove the uniqueness of 
solutions of the Ricci flow on complete non-compact manifolds by using 
the De Turck trick. However their proof is not correct because the 
crucial lemma Lemma 2.2 of \cite{CZ} is not correct. In Lemma 2.2 
of \cite{CZ} they 
claimed that they can construct a sequence of functions $\{\phi_a\}_{a\ge 1}$
which behaves like the distance function and have bounded covariant 
derivatives of all orders. They do this by smoothing the distance
function with the Riemannian convolution operator of R.E.~Green and H.~Wu 
(P.646--647 of \cite{GW1} and P.57 of \cite{GW2}). More precisely 
(\cite{GW1},\cite{GW2}) let 
$\psi:\Bbb{R}\to\Bbb{R}$ be a nonnegative $C^{\infty}$ function with support
in $[-1,1]$ which is constant in a neighborhood of $0$ and 
$$
\int_{v\in\Bbb{R}^n}\psi (|v|)=1.
$$
Then for any n-dimensional complete non-compact Riemannian manifold $(M, g)$ 
and continuous function $f:M\to\Bbb{R}$, the Riemannian convolution 
operator for $f$ is defined as
$$
f_{\3}(p)=\frac{1}{\3^n}\int_{v\in T_pM}f(exp_pv)\psi(|v|/\3)\,d\Omega_p
\quad\forall \3>0,p\in M
$$
where $d\Omega_p$ is the measure on $T_pM$ obtained from the Riemannian
metric on $M$ and $exp_p:T_pM\to M$ is the exponential map of $M$ at $p$. 
Hence
$$
f_{\3}(p)
=\frac{1}{\3^n}\int_{q\in M}f(q)\psi(|exp_p^{-1}(q)|/\3)\,dq\quad\forall
0<\3<\text{inj}(p),p\in M\tag 0.6
$$
As observed by R.E.~Green and H.~Wu (P.646--647 of \cite{GW1}) for the 
smoothness of $f_{\3}$ at $p\in M$, $2\3$ has to be less than the injectivity 
radius inj$(p)$ of $M$ at $p$. This is because one has to use the 
representation (0.6) for $f_{\3}$ in order to pass all the derivatives 
of $f_{\3}$ onto the smooth function
$$
\psi(|exp_p^{-1}(q)|/\3)
$$
in the integrand in (0.6). Thus one needs $exp_p:T_pM\to M$ to be 
a local diffeomorphism between the ball $B(0,2\3)$ in $T_pM$ and the
ball $B_g(p,2\3)\subset M$ for $0<2\3<\text{inj}(p)$, $p\in M$. Hence 
$$
|\nabla^kf_{\3}|(p)\approx C\3^{-k}\ge C\,\text{inj}(p)^{-k}\quad\forall
0<2\3<\text{inj}(p),p\in M, k\in\Bbb{Z}^+.\tag 0.7
$$
Let $p_0\in M$ be a fixed point of $M$ and suppose that $M$  also has bounded
curvature. Since the injectivity radius inj$(p)$ may 
decreases to $0$ as $dist (p_0,p)\to\infty$ \cite{CGT}, \cite{CLY}, by (0.7)
$|\nabla^kf_{\3}|(p)$ is not uniformly bounded on $M$ in general for any
$k\in\Bbb{Z}^+$. Thus the sequence of functions $\{\phi_a\}_{a\ge 1}$ 
constructed in \cite{CZ} can behave like
$$
|\nabla^k\phi_a|(p)\approx C\,\text{inj}(p)^{-k}
\quad\forall p\in M,k\in\Bbb{Z}^+
$$
and tends to infinity as $dist (p_0,p)\to\infty$. Hence Lemma 2.2 of 
\cite{CZ} is not correct. 

A simple example of manifold with bounded curvature and 
injectivity radii tending to zero as the point goes to infinity is as follows.
Let $f\in C^{\infty}(\R)$ be such that $f(x)>0$ on $\R$ and
$f(x)=|x|^{-2}$ for all $|x|\ge 1$. Consider the manifold
$M=S^1\times\R$ where $S^1$ is the circle with metric $g=f(x)d\theta^2+dx^2$.
Using the indices $1,2$ for the variable $x$ and $\theta$ respectively we
have
$$
\quad g_{1,1}(x,\theta)=1,\quad g_{2,2}(x,\theta)=f(x),\quad g_{2,1}(x,\theta)
=g_{1,2}(x,\theta)=0,\quad\forall\theta\in S^1,x\in\R.
$$
By direct computation,
$$
\Gamma_{11}^1=\Gamma_{11}^2=\Gamma_{21}^1=\Gamma_{12}^1=\Gamma_{22}^2=0,\quad
\Gamma_{22}^1=\frac{1}{x^3},\quad
\quad\Gamma_{21}^2=\Gamma_{12}^2=-\frac{1}{x}\quad\forall |x|\ge 1,
\theta\in S^1.
$$
Hence 
$$
R_{1212}=<(\nabla_1\nabla_2-\nabla_2\nabla_1)\frac{\1}{\1 x^1},\frac{\1}{\1 x^2}>
=\frac{2}{x^2}g_{22}=\frac{2}{x^4}\quad\forall |x|\ge 1,
\theta\in S^1.
$$
Thus
$$
|Rm|^2=2(g^{11}g^{22}R_{1212})^2=\frac{8}{x^4}\quad\forall |x|\ge 1,
\theta\in S^1.
$$
Hence 
$$
\|Rm\|_{L^{\infty}(M,g)}\le k_0<\infty
$$
for some constant $k_0>0$. On the other hand the injectivity radius 
$i(x,\theta)$ at a point $(x,\theta)\in M=S^1\times\R$ is less than or 
equal to the conjugate radius $\text{conj}\,(x,\theta)$ at $(x,\theta)$. 
Since $\text{conj}\, (x,\theta)=\pi f(x)=\pi/x^2\to 0$ as $|x|\to\infty$. 
Hence 
$$
i(x,\theta)\to 0\quad\text{ as }|x|\to\infty.
$$
Thus by this example and the above discussion Lemma 2.2 of \cite{CZ} wrong. 
Since Chen-Zhu's paper \cite{CZ} depends entirely on the validity of Lemma 2.2 
of \cite{CZ}. Hence Chen-Zhu's paper \cite{CZ} is wrong.  

Similarly according to the results and examples in \cite{CG}
there are many examples of sequences of manifolds with uniformly 
bounded curvature but with the corresponding injectivity radii converging
uniformly to zero. One can also read the survey article \cite{G} by 
J.D.E.~Grant on the injectivity radius estimate, the paper by 
J.~Cheeger, M.~Gromov and M.~Taylor \cite{CGT}, and (i) of Remark 1.7 
of \cite{AM} for various lower bound estimates on a manifold under 
various curvature conditions. Hence Lemma 2.2 of 
\cite{CZ} cannot be correct. 

Observe that the proof of \cite{CZ} uses the uniform boundedness 
property of the higher order covariant derivatives of the approximate 
distance function of Lemma 2.2 of \cite{CZ} in an essential way.
In this paper instead of Lemma 2.2 of \cite{CZ} I will use Corollary 1.5 
of this paper in the Deturck program to solve the uniqueness problem. 
Because of the absence of the uniform boundedness property of third and
higher order derivatives of the approximate distance function in 
Corollary 1.5, many theorems in 
this paper require new proofs different from that of \cite{CZ}. 

In the book \cite{SY} by R.~Schoen and S.T.~Yau a weaker result similar 
to Corollary 1.5 is proved in Theorem 4.2 of Chapter 1 using P.D.E. methods. 
However my proof of Corollary 1.5 is more elementary and requires only 
knowledge of the distance function on the manifold. 

Note that the manifold under consideration is non-compact and Lemma 3.5 
of \cite{H2} is applicable to the proof of the inequality on the last 
two lines of P.144 of \cite{CZ} only when $M$ is compact or the extremum 
of the norm of the covariant derivatives of the solution of Ricci harmonic 
map (0.2) can be attained in a compact set of $M$ independent of time.
Hence the proof of the uniform estimates over $M$ for the norms of the 
covariant derivatives of solutions of the harmonic map and the uniform 
lower bound estimate for the existence time of the solutions of the
approximate problems in Theorem 2.6 of \cite{CZ} on P.144--145 of \cite{CZ} 
is also not correct. Thus the proof of the existence of the Ricci harmonic
map in Theorem 2.1 of \cite{CZ} is not correct. 

The proof of Proposition 3.1 of \cite{CZ} has gaps and
the proof of Proposition 3.3 of \cite{CZ} which is crucial to the proof 
of uniqueness of solutions of Ricci flow is also not correct since the 
deduction of the last two inequalities on P.151 for the proof of 
Proposition 3.3 of \cite{CZ} assumed that one can interchange 
differentiation and taking limit as $\3\to 0$ which is also not true 
in general. 

In this paper we will give  a correct proof of the uniqueness 
of solutions of the Ricci flow on complete noncompact manifolds with bounded 
curvatures. We will use the De Turck approach to prove this result. We will 
prove the existence of solution of the Ricci harmonic flow on complete 
noncompact manifolds with bounded curvatures. 

The plan of the paper is as follows. In section 1 we will prove various 
estimates for the Hessian of the distance functions in both $(M,g)$ and 
the target manifold $(N,h)$. We will construct $C^2$ functions on $M$
with uniformly bounded first and second order covariant derivatives
which approximate the distance function of $(M,g(0))$. In section 2 we will 
construct solutions of (0.2) in bounded cylindrical domains with Dirichlet 
boundary condition and in $M\times (0,T_1)$ for some constant $T_1>0$.
We will prove the uniform estimates on the norm of the covariant derivatives 
of the solutions of the Ricci harmonic flow. In section 3 we will prove the 
uniqueness of the solutions of Ricci flow on complete noncompact manifolds 
with bounded curvatures.

We will let $(M,g(t))$, $0\le t\le T$, be a solution of the Ricci
flow on a $n$-dimensional complete non-compact manifold and 
$(N,h_{\alpha \beta})=(M,g(0))$ for the rest of the paper. 
We will assume that there exists a constant $k_0>0$ such that
$$
|\text{Rm}|\le k_0\quad\text{ on }M\times [0,T]\tag 0.8
$$
where Rm is the Riemannian curvature of $g(t)$ and $|\cdot|$ is the norm 
with respect to the metric $g(t)$. Note that by the results of W.X.~Shi
\cite{S1} for any $m\in\Bbb{Z}^+$ there exists a constant $c_m>0$ such 
that
$$
|\nabla^mR_{ijkl}|\le c_mt^{-\frac{m}{2}}.\tag 0.9
$$
For any $p,q\in M$, we let $\rho(p,q)$ be the distance between $p$ and $q$ 
with respect to $g(0)$. For any $k>0$, $p\in M$, let $B(p,k)
=\{q\in M:\rho (p,q)<k\}$. We 
will fix a point $p_0\in M$. For any $k>0$, let $B_k=B(p_0,k)$ and 
$Q_k^{T_1}=B(p_0,k)\times (0,T_1)$ for any $T_1>0$. For any bounded domain 
$\Omega\subset M$ we let 
$$
\1_p(\Omega\times (0,T_1))=\2{\Omega}\times\{0\}\cup\1\Omega\times [0,T_1)
$$  
be the parabolic boundary of $\Omega\times (0,T_1)$. For any open set
$O\subset M$, we let Vol$_{g(t)}(O)$ be the volume of $O$ with respect 
to the metric $g(t)$. For any $r>0$, let $V_{-k_0}(r)$ be the volume of 
a geodesic ball of radius $r$ in a space form of curvature $-k_0$.
For any point $x\in M$ we let Cut$(x)$ be the set of all cut points 
of $x$ with respect to the metrics $g(0)$.

Let $\nabla,\nabla^0,\4{\nabla},\hat{\nabla}$ be the 
covariant derivatives with respect the the metric $g(t)$,  $g(0)$, $h$, and 
$\hat{g}$ respectively. Let $\Gamma_{ij}^k(t)$, $\4{\Gamma}_{ij}^k(t)$, 
$R_{ijkl}$, $\4{R}_{ijkl}$, $R_{ij}$, $\4{R}_{ij}$, be Christoffel symbols, 
curvature tensors and Ricci tensors with respect to the metric $g(t)$ and 
$h$ respectively. Let $\Delta_t$ be the Laplace operator with respect to the
metric $g(t)$. 

We let $f:(M,g(0))\to (N,h)$ be a given diffeomorphism satisfying
$$
K_1=\sup_M|\nabla f|_{g(0),h}=\sup_M\biggl (g^{ij}(0)h_{\alpha\beta}
\frac{\1 f^{\alpha}}{\1 x^i}\frac{\1 f^{\beta}}{\1 x^j}\biggr )^{\frac{1}{2}}
<\infty\tag 0.10
$$
and
$$
K_2=\sup_M|{\nabla}^2f|_{g(0),h}<\infty.\tag 0.11
$$
for the rest of the paper. When there is no ambiguity we will drop the 
subscript and write $|\nabla f|$, $|{\nabla}^2f|$, for 
$|\nabla f|_{g(0),h}$ and $|{\nabla}^2f|_{g(0),h}$.

By the discussion in \cite{J} (cf. \cite{H3}) we can write the
derivative $\nabla F$ as
$$
\nabla F=\frac{\1 F^{\alpha}}{\1 x^i}dx^i\otimes\frac{\1}{\1 y^{\alpha}}
$$
and consider $\nabla F$ as a section of the bundle $T^{\ast}M\otimes 
F^{-1}TN$. The connection $\nabla$ on $(M,g(t))$ and the connection 
$\4{\nabla}$ on $(N,h)$ induce a natural connection $\nabla$ on  
$T^{\ast}M\otimes F^{-1}TN$ which can in turn be extended 
naturally to a connection on $T^{\ast}M\otimes T^{\ast}M\otimes F^{-1}TN$. 
By induction we get a natural connection on $T^{\ast}M^{\otimes p}
\otimes F^{-1}TN$ for any $p\in\Bbb{Z}^+$ from the connection $\nabla$ on 
$(M,g(t))$ and the connection $\4{\nabla}$ on $(N,h)$. 
More precisely for any (cf. \cite{J}, \cite{H3}, \cite{CZ})
$$
u=u_{i_1,i_2,\dots, i_{p-1}}^{\alpha}dx_{i_1}\otimes\cdots dx_{i_{p-1}}
\otimes\frac{\1}{\1 y^{\alpha}}\in T^{\ast}M^{\otimes {(p-1)}}\otimes F^{-1}TN
,p\ge 2,\tag 0.12
$$
$$
\nabla u=(\nabla_{i_p}u_{i_1,i_2,\dots,i_{p-1}}^{\alpha})\,
dx_{i_1}\otimes\cdots dx_{i_p}\otimes\frac{\1}{\1 y^{\alpha}}
\in T^{\ast}M^{\otimes p}\otimes F^{-1}TN
$$
where
$$
\nabla_{i_p}u_{i_1,i_2,\dots,i_{p-1}}^{\alpha}=\frac{\1}{\1 x^{i_p}}
u_{i_1,i_2,\dots, i_{p-1}}^{\alpha}
-\Gamma_{i_p,i_j}^mu_{i_1,\dots, i_{j-1},m,i_j,\dots, i_{p-1}}^{\alpha}
+\4{\Gamma}_{\beta\gamma}^{\alpha}\frac{\1 F^{\beta}}{\1 x^{i_p}}
u_{i_1,\dots, ,i_{p-1}}^{\gamma}.
$$
For the interchange of two covariant derivatives on 
$T^{\ast}M^{\otimes p}\otimes F^{-1}TN$ we have (cf. P.258 of \cite{H1} and 
P.133 of \cite{CZ})
$$
\nabla_i\nabla_ju_{i_1,\dots,i_p}^{\alpha}
-\nabla_j\nabla_iu_{i_1,\dots,i_p}^{\alpha}
=R_{iji_ml}g^{lm}u_{i_1,\dots, m,\dots,i_p}^{\alpha}
+\4{R}_{pqrs}\frac{\1 F^p}{\1 x^i}\frac{\1 F^q}{\1 x^j}h^{\alpha r}
u_{i_1,\dots,i_p}^s.
$$
We will equip $T^{\ast}M^{\otimes p}\otimes F^{-1}TN$ with the norm
$g^{\otimes p}\otimes h$ induced from $g$ and $h$. We write $|u|_{g(t),h}$
for the norm of $u\in T^{\ast}M^{\otimes p}\otimes F^{-1}TN$. When there
is no ambiguity we will drop the subscript and write $|u|$ instead of 
$|u|_{g(t),h}$. By abuse of notation we will also denote the operator
$g^{ij}\nabla_i\nabla_j$ on $u_{i_1,\dots,i_p}^{\alpha}$
by $\Delta_t$.
  
Similarly (cf. P.133 of \cite{CZ}) the time derivative $\1/\1 t$ can be 
extended naturally to a covariant time derivative $\nabla_t$ on 
$T^{\ast}M^{\otimes p}\otimes F^{-1}TN$ for any $p\ge 1$. For any 
$u\in T^{\ast}M^{\otimes {p}}\otimes F^{-1}TN$, $p\ge 1$ given by (0.12),
we define
$$
\nabla_tu=(\nabla_tu_{i_1,i_2,\dots, i_p}^{\alpha})\,
dx_{i_1}\otimes\cdots dx_{i_p}\otimes\frac{\1}{\1 y^{\alpha}}
$$
where
$$
\nabla_tu_{i_1,i_2,\dots, i_p}^{\alpha}=\frac{\1}{\1 t}
u_{i_1,i_2,\dots, i_p}^{\alpha}
+u_{i_1,\dots, ,i_p}^{\gamma}
\4{\Gamma}_{\beta\gamma}^{\alpha}\frac{\1 F}{\1 t}^{\beta}.
$$
Let $\eta\in C^{\infty}(\Bbb{R})$ , $0\le\eta\le 1$, be an even 
function such that $\eta (x)\equiv 1$ for $|x|\le 1/2$, $\eta (x)\equiv 0$ 
for $|x|\ge 1$, and $\eta'(x)\le 0$ for $x\ge 0$. Let 
$\phi=\eta^2$. Then $\phi\in C^{\infty}(\Bbb{R})$ is an even function, 
$0\le\phi\le 1$, $\phi (x)\equiv 1$ for $|x|\le 1/2$, $\phi (x)\equiv 0$ for 
$|x|\ge 1$, $\phi'(x)\le 0$ for $x\ge 0$, and
$$
\sup_{\Bbb{R}}({\phi'}^2/\phi)=4\sup_{\Bbb{R}}{\eta'}^2<\infty.\tag 0.13
$$ 

$$
\text{Section 1}
$$

In this section we will use a modification of the method of W.X.~Shi 
\cite{S2} to construct $C^2$ functions on $M$ with uniformly bounded first 
and second order covariant derivatives which approximate the distance 
function of $(M,g(0))$. We will obtain a sequence of regularization 
$\{\oa{h}\}$ for the metric $h$ each of which has a uniform lower bound on the 
injectivity radius on $N$. We first prove some estimates for the Hessian 
of the distance functions on $M$.

\proclaim{\bf Lemma 1.1}
Let $y_1,y_2\in M$ with $y_2\not\in\text{Cut}\,(y_1)$ and
$$
\rho(y_1,y_2)\le\pi/4\sqrt{k_0}\tag 1.1
$$
Let $\gamma$ be the unique minimal geodesic in $(M,g(0))$ from $y_1$ to 
$y_2$. Then for any unit vector $X\in T_{y_2}M$ perpendicular to 
$\1/\1\gamma$,
$$
Hess_{g(0)}(\rho)(X,X)\ge\frac{\pi}{4\rho(y_1,y_2)}\tag 1.2
$$
where $\rho=\rho(y_1,y_2)$. Hence for any unit vector $X\in 
T_{y_2}M$,
$$\left\{\aligned
&Hess_{g(0)}(\rho^2/2)(X,X)\ge\frac{\pi}{4}\\
&Hess_{g(0)}(\rho)(X,X)\ge 0.\endaligned\right.\tag 1.3
$$
\endproclaim
\demo{Proof}
Let $N_1$ be a space form of curvature $k_0$. Let $\rho_{N_1}(z)$ 
be the distance function on $N_1$ with respect to some fixed point 
$z_1\in N_1$. Suppose $z_2$ is a point on $N_1$ such that $\rho(y_1,y_2)
=\rho_{N_1}(z_2)$. Let $\rho=\rho(y_1,y_2)$ and let $\gamma_1$ be the 
minimal geodesics in $N_1$ from $z_1$ to $z_2$. Let $\xi\in T_{z_2}N_1$ 
be a unit vector which satisfies $<\xi,\1/\1\gamma_1>=0$. 
We extend $\xi$ to a vector field $X_1$ perpendicular to 
$\1/\1\gamma_1$ along $\gamma_1$ by parallel translation. Let
$$
f(s)=\frac{\sin \sqrt{k_0}s}{\sin \sqrt{k_0}\rho}.
$$
Then by (1.1) and an argument similar to that of \cite{SY}, 
$Y(s)=f(s)X_1(s)$ is the Jacobi field with $Y(\rho_0)=\xi$ and
$$\align
Hess_{g(0)}(\rho_{N_1})(\xi,\xi)(z)
=&\int_0^{\rho}\biggl (\biggl |\frac{df}{ds}\biggr |^2-k_0f(s)^2
\biggr )\,ds\\
=&\frac{k_0}{\sin^2\sqrt{k_0}\rho}\int_0^{\rho}(\cos^2\sqrt{k_0}s
-\sin^2\sqrt{k_0}s)\,ds\\
=&\frac{\sqrt{k_0}\sin 2\sqrt{k_0}\rho}
{2\sin^2\sqrt{k_0}\rho}\\
=&\frac{1}{\rho}\cdot
\frac{\sqrt{k_0}\rho}{\tan\sqrt{k_0}\rho}\\
\ge&\frac{1}{\rho}\cdot\frac{\pi/4}{\tan(\pi/4)}\\
\ge&\frac{\pi}{4\rho(y_1,y_2)}.\tag 1.4
\endalign
$$ 
Hence by (0.8), (1.4), and the Hessian comparsion theorem \cite{SY}, 
for any unit vector $X\in T_{y_2}M$ satisfying $<X,\1/\1\gamma>=0$, we have 
$$\align
&Hess_{g(0)}(\rho)(X,X)\ge Hess_{g(0)}(\rho_{N_1})(\xi,\xi)
\ge\frac{\pi}{4\rho(y_1,y_2)}\\
\Rightarrow\quad&Hess_{g(0)}(\rho^2/2)(X,X)\ge\frac{\pi}{4}\tag 1.5
\endalign
$$
and (1.2) follows. By direct computation, for any $X\in T_{y_2}M$ satisfying 
$<X,\1/\1\gamma>=0$,
$$\left\{\aligned
&Hess_{g(0)}(\rho)(\1/\1\gamma,\1/\1\gamma)
=Hess_{g(0)}(\rho)(X,\1/\1\gamma)=Hess_{g(0)}(\rho)(\1/\1\gamma,X)
=0\\
&Hess_{g(0)}(\rho^2/2)(X,\1/\1\gamma)=Hess_{g(0)}(\rho^2/2)
(\1/\1\gamma,X)=0\\
&Hess_{g(0)}(\rho^2/2)(\1/\1\gamma,\1/\1\gamma)=1.
\endaligned\right.\tag 1.6
$$
For any $X\in T_{y_2}M$,
$$
X=X_1+\lambda\frac{\1}{\1\gamma}
$$
for some constant $\lambda$ and $X_1\in T_{y_2}M$ perpendicular to 
$\1/\1\gamma$. Then by (1.5) and (1.6), $\forall X\in 
T_{y_2}M$,
$$\align
Hess_{g(0)}(\rho)(X,X)=&Hess_{g(0)}(\rho)(X_1,X_1)+2\lambda 
Hess_{g(0)}(\rho)(X,\frac{\1}{\1\gamma})\\
&\qquad+\lambda^2Hess_{g(0)}(\rho)(\frac{\1}{\1\gamma},\frac{\1}{\1\gamma})\\
\ge&0
\endalign
$$
and
$$\align
Hess_{g(0)}(\rho^2/2)(X,X)=&Hess_{g(0)}(\rho^2/2)(X_1,X_1)+2\lambda 
Hess_{g(0)}(\rho^2/2)(X,\frac{\1}{\1\gamma})\\
&+\lambda^2Hess_{g(0)}(\rho^2/2)(\frac{\1}{\1\gamma},\frac{\1}{\1\gamma})\\
\ge&\frac{\pi}{4}|X_1|_{g(0)}^2+\lambda^2\\
\ge&\frac{\pi}{4}
\endalign
$$
and (1.3) follows. 
\enddemo

\proclaim{\bf Lemma 1.2}(cf. P.225-226 of \cite{S1})
There exist constants $c_2>c_1>0$ such that
$$\left\{\aligned
&c_1g_{ij}(x,0)\le g_{ij}(x,t)\le c_2g_{ij}(x,0)\quad\forall 
0\le t\le T\\
&c_1g^{ij}(x,0)\le g^{ij}(x,t)\le c_2g^{ij}(x,0)\quad\forall 
0\le t\le T.\endaligned\right.\tag 1.7
$$
\endproclaim

\proclaim{\bf Lemma 1.3}
Let $x_0,x\in B_k$ with $x\not\in\text{Cut}\,(x_0)$ and let 
$\rho(x)=\rho(x_0,x)$ satisfy (1.1). Then exists a constant $C_1>0$ such 
that
$$
-C_1\sqrt{t}\le\Delta_t\rho\le C_1(1+\frac{1}{\rho})
\quad\forall 0\le t\le T.\tag 1.8
$$
\endproclaim
\demo{Proof}
Note that
$$
\Delta_t\rho=g^{ij}(t)\nabla_i\nabla_j\rho
=g^{ij}(t)\biggl (\nabla_i^0\nabla_j^0\rho-(\Gamma_{ij}^k(t)-\Gamma_{ij}^k(0))
\frac{\1\rho}{\1 x^k}\biggr ).\tag 1.9
$$
We choose a normal coordinate system $\{\1/\1 x^i\}$ with respect to the
metric $g_{ij}(t)$ at $x$. Then by (0.1), (0.8), (0.9), Lemma 1.2, and an 
argument similar to the proof of (2.16) and (2.17) of \cite{CZ},
$$
\biggl |\frac{\1}{\1 t}\Gamma_{ij}^k(t)\biggr |\le\frac{C}{\sqrt{t}}
\quad\Rightarrow\quad
|\Gamma_{ij}^k(t)-\Gamma_{ij}^k(0)|\le C\sqrt{t}
\quad\forall 0\le t\le T.\tag 1.10
$$
By the Hessian comparison theorem \cite{SY}, P.309-310 of \cite{S2},
and Lemma 1.2 (cf. P.136 of \cite{CZ}),
$$\align
&\nabla_i^0\nabla_j^0\rho\le\frac{1+\sqrt{k_0}\rho}{\rho}g_{ij}(x,0)\\
\Rightarrow\quad&g^{ij}(x,t)\nabla_i^0\nabla_j^0\rho
\le c_2\frac{1+\sqrt{k_0}\rho}{\rho}g^{ij}(x,0)g_{ij}(x,0)
\le nc_2\biggl (\frac{1+\sqrt{k_0}\rho}{\rho}\biggr ).\tag 1.11
\endalign
$$
By Lemma 1.1 (1.3) holds. Hence
$$
g^{ij}(x,t)\nabla_i^0\nabla_j^0\rho=\nabla_i^0\nabla_i^0\rho\ge 0.\tag 1.12
$$
By (1.9), (1.10), (1.11) and (1.12) we get (1.8) and the lemma follows.
\enddemo

Let 
$$
k_1=\pi/4\sqrt{k_0}
$$ 
and
$$
\2{\rho}(y)=\2{\rho}(p_0,y)=\frac{\int_M\rho(p_0,z)
\eta (\rho(y,z)/k_1)\,dz}{\int_M\eta (\rho(y,z)/k_1)\,dz}.
$$
Then 
$$
\rho(p_0,y)-k_1\le\2{\rho}(y)\le\rho(p_0,y)+k_1
\quad\forall y\in M.\tag 1.13
$$

\proclaim{\bf Lemma 1.4}
$\2{\rho}\in C^2(M)$ and there exists a constant $C_1>$ such that
$$\left\{\aligned
&|\4{\nabla}\2{\rho}(y)|_h\le C_1\quad\forall y\in M\\
&|\4{\nabla}^2\2{\rho}(y)|_h\le C_1\quad\forall y\in M
\endaligned\right.
$$
\endproclaim
\demo{Proof}
Since 
$$
z\in M\setminus\text{Cut}\,(y)\quad\Leftrightarrow\quad
y\in M\setminus\text{Cut}\,(z),\tag 1.14
$$ 
$$
|\4{\nabla}\rho(y,z)|_h\le 1\quad\forall z\in M\setminus\text{Cut}\,(y).
\tag 1.15
$$
As Cut$(y)$ has measure zero, by (1.15) and the Lebesgue dominated 
convergence theorem,
$$
\4{\nabla}_{\beta}\2{\rho}(y)=k_1^{-1}
\frac{\int_M(\rho(p_0,z)-\2{\rho}(y))\eta'\4{\nabla}_{\beta}
\rho(y,z)\,dz}{\int_M\eta (\rho(y,z)/k_1)\,dz}
\tag 1.16
$$
where $\eta'=\eta'(\rho(y,z)/k_1)$. By (1.13), (1.15), (1.16) and the 
volume comparison theorem \cite{SY}, \cite{Ch}, there exist constants 
$C_1>0$ and $C_1'>0$ such that
$$
|\4{\nabla}\2{\rho}(y)|_h\le (C_1'/k_1)
\frac{\text{Vol}_h(B_h(p_0,k_1))}{\text{Vol}_h(B_h(p_0,k_1/2))}\cdot
\sup_{z\in B_h(y,k_1)}|\rho(p_0,z)-\2{\rho}(y)|
\le 2C_1'\frac{V_{-k_0}(k_1)}{V_{-k_0}(k_1/2)}\le C_1\tag 1.17
$$
for any $y\in M$. By the Hessian comparison theorem \cite{SY}, Lemma 1.1 
(cf. P.309--311 of \cite{S2}), and (1.14),
$$
\4{\nabla}_{\alpha}\4{\nabla}_{\beta}\rho(y,z)
\le\frac{1+\sqrt{k_0}\rho}{\rho}h_{\alpha\beta}
\quad\forall z\in M\setminus\text{Cut}(y)
\tag 1.18
$$
and
$$\align
&Hess_{g(0)}(\rho(y,z))(\frac{\1}{\1 y^{\alpha}}+\frac{\1}{\1 y^{\beta}},
\frac{\1}{\1 y^{\alpha}}+\frac{\1}{\1 y^{\beta}})\ge 0\\
\Rightarrow\quad&Hess_{g(0)}(\rho(y,z))(\frac{\1}{\1 y^{\alpha}},
\frac{\1}{\1 y^{\beta}})\\
\ge&-\frac{1}{2}\biggl (Hess_{g(0)}(\rho(y,z))
(\frac{\1}{\1 y^{\alpha}},\frac{\1}{\1 y^{\alpha}})
+Hess_{g(0)}(\rho(y,z))(\frac{\1}{\1 y^{\beta}},\frac{\1}{\1 y^{\beta}})
\biggr )\tag 1.19
\endalign
$$
holds for any $z\in M\setminus\text{Cut}\,(y)$. We now choose a normal 
co-ordinate system $\{\1/\1y^{\alpha}\}$ at $y$. Then by (1.15), (1.18),
and (1.19),
$$\align
&|\4{\nabla}_{\alpha}\4{\nabla}_{\beta}\rho|_h
\le n^2\biggl(\frac{1+\sqrt{k_0}\rho}{\rho}\biggr)
\quad\forall z\in M\setminus\text{Cut}(y)\tag 1.20\\
\Rightarrow\quad&|\rho(p_0,z)(\eta''\4{\nabla}_{\alpha}
\rho\4{\nabla}_{\beta}
\rho+k_1\eta'\4{\nabla}_{\alpha}\4{\nabla}_{\beta}\rho)|_h
\le C_2(\rho(p_0,y)+1)\chi_{B_h(y,k_1)}\tag 1.21
\endalign
$$
holds for any $z\in M\setminus\text{Cut}(y)$ where $\eta''
=\eta''(\rho(y,z)/k_1)$ and $\chi_{B_h(y,k_1)}$ is the 
characteristic function of the set $B_h(y,k_1)$. Hence by (1.16), (1.21), 
and the Lebesgue dominated convergence theorem,
$$\align
\4{\nabla}_{\alpha}\4{\nabla}_{\beta}\2{\rho}(y)=&k_1^{-2}\biggl\{
\frac{\int_M\rho(p_0,z)(\eta''\4{\nabla}_{\alpha}\rho\4{\nabla}_{\beta}
\rho+k_1\eta'\4{\nabla}_{\alpha}\4{\nabla}_{\beta}\rho)\,dz}
{\int_M\eta (\rho(y,z)/k_1)\,dz}\\
&-\frac{\int_M\rho(p_0,z)\eta'\4{\nabla}_{\beta}
\rho\,dz}{\int_M\eta (\rho(y,z)/k_1)\,dz}\cdot
\frac{\int_M\eta'\4{\nabla}_{\alpha}\rho\,dz}
{\int_M\eta (\rho(y,z)/k_1)\,dz}\\
&-k_1\4{\nabla}_{\alpha}\2{\rho}(y)
\frac{\int_M\eta'\4{\nabla}_{\beta}\rho\,dz}
{\int_M\eta (\rho(y,z)/k_1)\,dz}\\
&-\2{\rho}(y)\frac{\int_M(\eta''\4{\nabla}_{\alpha}\rho\4{\nabla}_{\beta}
\rho+k_1\eta'\4{\nabla}_{\alpha}\4{\nabla}_{\beta}\rho)\,dz}
{\int_M\eta (\rho(y,z)/k_1)\,dz}\\
&-\2{\rho}(y)\frac{\int_M\eta'\4{\nabla}_{\alpha}
\rho\,dz}{\int_M\eta (\rho(y,z)/k_1)\,dz}\cdot
\frac{\int_M\eta'\4{\nabla}_{\beta}\rho\,dz}
{\int_M\eta (\rho(y,z)/k_1)\,dz}\biggr\}\\
=&k_1^{-2}\biggl\{
\frac{\int_M(\rho(p_0,z)-\2{\rho}(y))
(\eta''\4{\nabla}_{\alpha}\rho\4{\nabla}_{\beta}
\rho+k_1\eta'\4{\nabla}_{\alpha}\4{\nabla}_{\beta}\rho)\,dz}
{\int_M\eta (\rho(y,z)/k_1)\,dz}\\
&-\frac{\int_M(\rho(p_0,z)-\2{\rho}(y))\eta'\4{\nabla}_{\beta}
\rho\,dz}{\int_M\eta (\rho(y,z)/k_1)\,dz}\cdot
\frac{\int_M\eta'\4{\nabla}_{\alpha}\rho\,dz}
{\int_M\eta (\rho(y,z)/k_1)\,dz}\\
&-k_1\4{\nabla}_{\alpha}\2{\rho}(y)
\frac{\int_M\eta'\4{\nabla}_{\beta}\rho\,dz}
{\int_M\eta (\rho(y,z)/k_1)\,dz}\biggr\}\\
\tag 1.22
\endalign
$$
where $\rho=\rho(y,z)$. By (1.13), (1.17), (1.20), (1.21), (1.22), 
and the volume comparison theorem, 
$$\align
&|\4{\nabla}_{\alpha}\4{\nabla}_{\beta}\2{\rho}(y)|_h\\
\le&C_2'\biggl\{\frac{\text{Vol}_h(B_h(p_0,k_1))}{\text{Vol}_h(B_h(p_0,k_1/2))}
+\biggl(\frac{\text{Vol}_h(B_h(p_0,k_1))}{\text{Vol}_h(B_h(p_0,k_1/2))}
\biggr)^2\biggr\}(1+\sup_{z\in B_h(y,k_1)}|\rho(p_0,z)-\2{\rho}(y)|)\\
\le&(2k_1+1)C_2'\biggl\{\frac{V_{-k_0}(k_1)}{V_{-k_0}(k_1/2)}+\biggl(
\frac{V_{-k_0}(k_1)}{V_{-k_0}(k_1/2)}\biggr)^2\biggr\}\\
\le&C_1'\quad\forall y\in M\tag 1.23
\endalign
$$
for some cosntants $C_1'>0$ and $C_2'>0$. By (1.16), (1.21) and (1.22), 
$\2{\rho}\in C^2(M)$ and the lemma follows.
\enddemo

For any $a\ge 1$, let
$$
\2{\rho}_a(y)=\2{\rho}(y)(1-\eta (\2{\rho}(y)/a)).\tag 1.24
$$
Then
$$
\2{\rho}_a(y)=0\quad\forall\2{\rho}(y)\le a/2.\tag 1.25
$$
Since $\2{\rho}_a(y)=\2{\rho}(y)$ for any $\2{\rho}(y)\ge a$, by (1.13),
$$
\rho(p_0,y)-k_1\le\2{\rho}_a(y)\le\rho(p_0,y)+k_1
\quad\forall \2{\rho}(y)\ge a.\tag 1.26
$$
Since
$$
\4{\nabla}_{\beta}\2{\rho}_a=(1-\eta(\2{\rho}/a))\4{\nabla}_{\beta}\2{\rho}
-\2{\rho}\eta'\cdot(\4{\nabla}_{\beta}\2{\rho}/a)
$$
and
$$
\4{\nabla}_{\alpha}\4{\nabla}_{\beta}\2{\rho}_a=(1-\eta(\2{\rho}/a))
\4{\nabla}_{\alpha}\4{\nabla}_{\beta}\2{\rho}
-2(\eta'/a)\4{\nabla}_{\alpha}\2{\rho}\4{\nabla}_{\beta}\2{\rho}
-\2{\rho}(\eta'/a)\4{\nabla}_{\alpha}\4{\nabla}_{\beta}\2{\rho}
-\2{\rho}(\eta''/a^2)\4{\nabla}_{\alpha}\2{\rho}\4{\nabla}_{\beta}\2{\rho},
$$
by Lemma 1.4 we have the following result.

\proclaim{\bf Corollary 1.5}
$\2{\rho}_a\in C^2(M)$ for any $a\ge 1$ and there exists a constant
$C_2>0$ such that
$$\left\{\aligned
&|\4{\nabla}\2{\rho}_a(y)|_h\le C_2\quad\forall y\in M,
a\ge 1\\
&|\4{\nabla}^2\2{\rho}_a(y)|_h\le C_2\quad\forall y\in M,
a\ge 1.\endaligned\right.
$$
\endproclaim

For any $a\ge 1$, let 
$$
\psi_a(y)=4\sqrt{k_0}\2{\rho}_a(y)
$$
and 
$$\left\{\aligned
&\oa{h}_{\alpha\beta}=e^{\psi_a}h_{\alpha\beta}\\
&\oa{g}_{ij}=e^{\psi_a}g_{ij}.
\endaligned\right.
$$
Note that by the results of \cite{CGT} and \cite{CLY}, there exists a 
constant $\delta_0>0$ depending on $k_0$ and the injectivity radius of 
$(N,h)$ at $p_0$ such that 
$$
inj_h(y)\ge\delta_0\,e^{-\sqrt{k_0}\rho(p_0,y)}\quad\forall y\in N
\tag 1.27
$$ 
where $inj_h(y)$ is the injectivity radius of $y$ in $(N,h)$.
Now by (1.13) and (1.26),
$$\align
&\rho(p_0,y)-k_1\le\2{\rho}_a(y)\le\rho(p_0,y)+k_1
\quad\forall\rho(p_0,y)\ge a+k_1,a\ge 1\\
\Rightarrow\quad&\2{\rho}_a(y)\ge 2\rho(p_0,y)/3\quad\forall
\rho(p_0,y)\ge 4a/3,a\ge\max (1,3k_1)\tag 1.28
\endalign
$$
where $k_1=\pi/4\sqrt{k_0}$. By (1.27), (1.28), and an argument similar 
to the proof on P.125 of \cite{CZ}, for any $a\ge\max (1,3\pi/4
\sqrt{k_0})$,
$$
i_a=\text{inj}\,(N,h^a)>0.\tag 1.29
$$
We will now let $\overset{a}\to{\text{Rm}}$, $\overset{a}\to{\4{\text{Rm}}}$, 
$\overset{a}\to{R}_{ijkl}$, $\overset{a}\to{\4{R}}_{\alpha\beta\gamma
\delta}$, etc. be the Riemann curvature, Riemannian curvature tensor, etc. 
of $\oa{g}_{ij}$ and $\oa{h}_{\alpha\beta}$ respectively. 

\proclaim{\bf Lemma 1.6}
There exists a constant $C_3>0$ such that
$$
|\overset{a}\to{\4{\text{Rm}}}|_{h^a}\le C_3e^{-\psi_a}\quad\text{ in }N
\quad\forall a\ge\max(1,3\pi/4\sqrt{k_0}).
$$
\endproclaim
\demo{Proof}
Let $a\ge\max(1,3\pi/4\sqrt{k_0})$. By direct computation
(cf. (2.9) of \cite{CZ} and (13) of P.299 of \cite{S1}),
$$\align
\overset{a}\to{\4{R}}_{\alpha\beta\gamma\delta}
=&e^{\psi_a}\4{R}_{\alpha\beta\gamma\delta}
+\frac{e^{\psi_a}}{4}\{|\4{\nabla}\psi_a|^2
(h_{\alpha\delta}h_{\beta\gamma}-h_{\alpha\gamma}h_{\beta\delta})
+(2\4{\nabla}_{\alpha}\4{\nabla}_{\delta}\psi_a-\4{\nabla}_{\alpha}\psi_a
\4{\nabla}_{\delta}\psi_a)h_{\beta\gamma}\\
&\qquad +(2\4{\nabla}_{\beta}\4{\nabla}_{\gamma}\psi_a
-\4{\nabla}_{\beta}\psi_a\4{\nabla}_{\gamma}\psi_a)h_{\alpha\delta}
-(2\4{\nabla}_{\beta}\4{\nabla}_{\delta}\psi_a-\4{\nabla}_{\beta}\psi_a
\4{\nabla}_{\delta}\psi_a)h_{\alpha\gamma}\\
&\qquad -(2\4{\nabla}_{\alpha}\4{\nabla}_{\gamma}\psi_a
-\4{\nabla}_{\alpha}\psi_a\4{\nabla}_{\gamma}\psi_a)h_{\beta\delta}\}.
\tag 1.30
\endalign
$$
Hence by Corollary 1.5 and (1.30),
$$\align
|\overset{a}\to{\4{\text{Rm}}}|_{h^a}\le&Ce^{-\psi_a}(|\4{Rm}|_{h}
+|\4{\nabla}^2\psi_a|_h+|\4{\nabla}\psi_a|_h^2)\\
=&Ce^{-\psi_a}(|\4{Rm}|_{h}+4|\4{\nabla}^2\2{\rho}_a|_h
+16|\4{\nabla}\2{\rho}_a|_h^2)\\
\le&C_3'e^{-\psi_a}(|\4{Rm}|_{h}+1)
\endalign
$$
for some constants $C>0, C_3'>0$, and the lemma follows.
\enddemo

By a similar argument as the proof of Lemma 1.6 we have the following result.

\proclaim{\bf Corollary 1.7}
There exists a constant $C_4>0$ such that
$$
|\overset{a}\to{\text{Rm}}|_{g^a}\le C_4e^{-\psi_a}\quad\text{ in }M\times
(0,T)\quad\forall a\ge\max(1,3\pi/4\sqrt{k_0}).
$$
\endproclaim

$$
\text{Section 2}
$$

In this section we will construct solutions $\overset{a}\to F:
(M,\ov{a}{g})\to (N,\ov{a}{h})$ of 
$$
\frac{\1\oa{F}}{\1 t}=\oa{\Delta}_{\oa{g}(t),\oa{h}}\oa{F}\tag 2.1
$$
in bounded cylindrical domains with Dirichlet boundary condition where
$$
\oa{\Delta}_{\oa{g}(t),\oa{h}}F=\Delta_{\oa{g}(t)}F+\oa{g}^{ij}(x,t)
\oa{\4{\Gamma}}^{\alpha}_{\beta ,\gamma}(\oa{F}(x,t))\frac{\1\oa{F}^\beta}
{\1 x^i}\frac{\1\oa{F}^\gamma}{\1 x^j}
$$
and 
$$
\Delta_{\oa{g}(t)}F^{\alpha}=\oa{g}_{ij}\oa{\nabla}_i\oa{\nabla}_jF^{\alpha}.
$$
We will construct solution for the approximate problem 
$$\left\{\aligned
&\frac{\1\oa{F}}{\1 t}=\oa{\Delta}_{\oa{g}(t),\oa{h}}F
\quad\text{ in }M\times (0, T_1)\\
&F(x,0)=f(x)\quad\text{ in }M
\endaligned\right.\tag 2.2
$$
of (0.2) in $M\times (0, T_1)$.
We will prove that the solutions of (2.2) has a subsequence that
converges to a solution of (0.2) as $a\to\infty$. We will now assume
that $a\ge\max(1,3\pi/4\sqrt{k_0})$ for the rest of the paper. Note that
as before the Levi-Civita connections $\oa{\nabla}$ and $\oa{\4{\nabla}}$
on $(M,\oa{g})$ and $(N,\oa{h})$ respectively induce a natural connection 
$\oa{\nabla}$ on $T^{\ast}M^{\otimes {p}}\otimes F^{-1}TN$ for any $p\ge 1$.
We will write $\oa{\Delta}_t$ for $\Delta_{\oa{g}(t)}
=\oa{g}_{ij}\oa{\nabla}_i\oa{\nabla}_j$. 

Note that 
$$
\sup_M|\oa{\nabla}f|_{\oa{g}(0),\oa{h}}=\sup_M|\nabla f|_{g(0),h}=K_1
\quad\forall a\ge 1.\tag 2.3
$$
By direct computation,
$$
\oa{\Gamma}_{ij}^k=\Gamma_{ij}^k+\frac{1}{2}\{\delta_i^k
\nabla_j\psi_a+\delta_j^k\nabla_i\psi_a-g^{kl}g_{ij}\nabla_l\psi_a\}.\tag 2.4
$$
Hence
$$\align
\oa{\nabla}_i\oa{\nabla}_jf^{\alpha}=&\nabla_i\nabla_jf^{\alpha}
+(\Gamma_{ij}^k-\oa{\Gamma}_{ij}^k)\nabla_kf^{\alpha}
+(\oa{\4{\Gamma}}_{\beta\gamma}^{\alpha}-\Gamma_{\beta\gamma}^{\alpha})
\nabla_if^{\beta}\nabla_jf^{\gamma}\\
=&\nabla_i\nabla_jf^{\alpha}-\frac{1}{2}\{\delta_i^k
\nabla_j\psi_a+\delta_j^k\nabla_i\psi_a-g^{kl}g_{ij}\nabla_l\psi_a\}
\nabla_kf^{\alpha}\\
&\qquad +\frac{1}{2}\{\delta_{\beta}^{\alpha}
\nabla_{\gamma}\psi_a+\delta_{\gamma}^{\alpha}\nabla_{\beta}\psi_a
-h^{\alpha\delta}h_{\beta\gamma}\nabla_{\delta}\psi_a\}
\nabla_if^{\beta}\nabla_jf^{\gamma}\\
=&\nabla_i\nabla_jf^{\alpha}-\frac{1}{2}\{\nabla_if^{\alpha}
\nabla_j\psi_a+\nabla_jf^{\alpha}\nabla_i\psi_a-g^{kl}g_{ij}\nabla_kf^{\alpha}
\nabla_l\psi_a\}\\
&\qquad +\frac{1}{2}\{\nabla_if^{\alpha}\nabla_jf^{\gamma}
\nabla_{\gamma}\psi_a+\nabla_jf^{\alpha}\nabla_if^{\beta}\nabla_{\beta}\psi_a
-h^{\alpha\delta}h_{\beta\gamma}\nabla_if^{\beta}\nabla_jf^{\gamma}
\nabla_{\delta}\psi_a\}\tag 2.5
\endalign
$$
By (0.10), (0.11), (2.5) and Corollary 1.5, there exists a constant $C>0$ 
such that
$$
|\oa{\nabla}^2f(x)|_{\oa{g}(0),\oa{h}}\le Ce^{-\frac{1}{2}\psi_a(x)}
\le C\quad\forall x\in M,a\ge 1.\tag 2.6
$$
By Lemma 1.6, Lemma 1.1 for $(N,\oa{h})$, and an argument similar 
to the proof of Lemma 2.8 of \cite{CZ} we have

\proclaim{\bf Lemma 2.1}
For any $a\ge\max(1,3\pi/4\sqrt{k_0})$, there exist constants $0<T_2\le T$
and $C_5>0$ depending on $n$, $a$, and $k_0$ such that for any $k>0$, 
$0<T_1\le T$, and solution $\oa{F}_k:(B_k,g(t))\to (N,h_{\alpha\beta})$, 
$\oa{F}_k\in C^{2+\frac{1}{2},1+\frac{1}{2}}(\2{B}_k\times [0,T_1))$, of
$$\left\{\aligned
&\frac{\1\oa{F}_k}{\1 t}=\oa{\Delta}_{\oa{g}(t),\oa{h}}
\oa{F}_k\quad\text{ in }Q_k^{T_1}\\
&\oa{F}_k(x,t)=f(x)\quad\text{ on }\1 B_k\times (0,T_1)\\
&\oa{F}_k(x,0)=f(x)\quad\text{ in }B_k,
\endaligned\right.\tag 2.7
$$
we have
$$
\rho(f(x),\oa{F}_k(x,t))\le C_5\sqrt{t}\tag 2.8
$$ 
in $Q_k^{T_2'}$ where $T_2'=\min (T_2,T_1)$.
\endproclaim

By Lemma 1.6, Lemma 1.7, (2.3), (2.6), an argument similar to the proof of 
Theorem 7.1 of Chapter VII of \cite{LSU} (cf. P.245--246 of \cite{S1}) 
and an argument similar to the proof of Lemma 2.9 of \cite{CZ} but with 
Lemma 2.1 replacing Lemma 2.8 in the proof there we have

\proclaim{\bf Lemma 2.2}
Let $a\ge\max(1,3\pi/4\sqrt{k_0})$ and let $T_2$ be as given in Lemma 2.1. 
Then there exist constants $0<T_3\le T_2$ depending on $n$, $a$, $k_0$, 
and $f$ such that for any $k>0$ there exists a solution $\oa{F}_k\in
C^{2+\frac{1}{2},1+\frac{1}{2}}(\2{B}_k\times [0,T_3])$ of (2.7) in 
$Q_k^{T_3}$.
\endproclaim

\proclaim{\bf Remark} 
Note that by Corollary 1.5 and an argument similar to the proof of 
Theorem 7.1 of Chapter VII of \cite{LSU} for any solution $\oa{F}_k$ 
of (2.7) in $Q_k^{T_3}$ given by Lemma 2.2, we have $\oa{F}_k\in
C^{3,1+\frac{1}{2}}(\2{B}_k\times (0,T_3])$.
\endproclaim

\proclaim{\bf Lemma 2.3}(Section 6 of \cite{H3})
Let $\oa{F}$ be a solution of (2.1) in $Q_k^{T'}$ for some $k>0$ and 
$0<T'\le T$. Then $\oa{F}$ satisfies
$$
\frac{\1}{\1 t}|\oa{\nabla}\oa{F}|^2=\oa{\Delta}_t|\oa{\nabla}\oa{F}|^2\
-2|\oa{\nabla}^2\oa{F}|^2+2\oa{\4{Rm}}
(\oa{\nabla}\oa{F},\oa{\nabla}\oa{F},\oa{\nabla}\oa{F},
\oa{\nabla}\oa{F})\quad\text{ in }Q_k^{T'}
$$
where
$$
\oa{\4{Rm}}
(\oa{\nabla}\oa{F},\oa{\nabla}\oa{F},\oa{\nabla}\oa{F},
\oa{\nabla}\oa{F})=\oa{g}^{ik}\oa{g}^{jl}\oa{\4{R}}_{pqrs}
\oa{\nabla}_i\oa{F}^p\oa{\nabla}_j\oa{F}^q\oa{\nabla}_k\oa{F}^r
\oa{\nabla}_l\oa{F}^s.
$$
\endproclaim

We next observe that by Lemma 1.6 and Lemma 2.3,
$$
\frac{\1}{\1 t}|\oa{\nabla}\oa{F}|^2\le\oa{\Delta}_t|\oa{\nabla}\oa{F}|^2
-2|\oa{\nabla}^2\oa{F}|^2+C_3|\oa{\nabla}F|^4\tag 2.9
$$
in $Q_k^{T'}$ for any solution $\oa{F}$ of (2.1) in $Q_k^{T'}$ and any 
$a\ge\max (1,3\pi/4\sqrt{k_0})$ where $C_3>0$ is the constant given by 
Lemma 1.6. Then by (2.3), (2.9), Lemma 2.1 and an argument 
similar to the proof of Lemma 2.11 of \cite{CZ} but with (2.9) replacing 
(2.30) in the proof there we have

\proclaim{\bf Lemma 2.4}
Let $T_3$ be as in Lemma 2.2. Let $k_0>0$ be given by (0.8).
Then there exist constants $0<\delta_1<1$, $C_6>0$, and $0<T_4\le T_3$ 
depending on $n$, $a$, $k_0$, and $K_1$ such that for any solution 
$\oa{F}_k$ of (2.7) given by Lemma 2.2 and any $B_{\oa{g}(0)}(x_0,\delta_1)
\subset B_k$ we have 
$$
|\oa{\nabla}\oa{F}_k|(x,t)\le C_6\tag 2.10
$$
holds in $B_{\oa{g}(0)}(x_0,3\delta_1/4)\times [0,T_4]$ for any $k>0$ and
$a\ge\max (1,3\pi/4\sqrt{k_0})$.
\endproclaim

\proclaim{\bf Theorem 2.5}
Let $a\ge\max (1,3\pi/4\sqrt{k_0})$ and let $T_4$, $C_6$, be given 
by Lemma 2.4. Then there exists a solution $\oa{F}\in C^{2+\frac{1}{2},
1+\frac{1}{2}}(M\times [0,T_4])\cap C^{3,1+\frac{1}{2}}(M\times 
(0,T_4])$ of (2.2) in $M\times (0,T_4)$ which satisfies
$$
|\oa{\nabla}\oa{F}|(x,t)\le C_6\quad\text{ on }M\times [0,T_4].\tag 2.11
$$
\endproclaim
\demo{Proof}
By Lemma 2.2 for any $k\in\Bbb{Z}^+$ there exists a solution $\oa{F}_k\in 
C^{2+\frac{1}{2},1+\frac{1}{2}}(\2{B}_k\times [0,T_4])$ of (2.7) in 
$Q_k^{T_4}$. By Lemma 2.4 (2.10) holds. Hence for any $k'\in\Bbb{Z}^+$ 
the sequence $\{|\oa{\nabla}\oa{F}_k|\}_{k=k'}^{\infty}$ is uniformly 
bounded on any compact subset of $\2{B_{\oa{g}(0)}(p_0,k'-1)}\times [0,T_4]$. 
Then by Corollary 1.5 and an argument similar 
to the proof of Theorem 7.1 of Chapter VII of \cite{LSU} (cf. P.245--246 of 
\cite{S1}) $\{\oa{F}_k\}_{k=k'}^{\infty}$ is uniformly bounded in 
$C^{2+\frac{1}{2},1+\frac{1}{2}}(\2{B_{\oa{g}(0)}(p_0,k'-1)}\times 
[0,T_4])$.  

By the Ascoli Theorem and a diagonalization argument the sequence 
$\{\oa{F}_k\}_{k=1}^{\infty}$ has a convergent subsequence
which we may assume without loss of generality to be the sequence itself
such that $\oa{F}_k$ converges uniformly in $C^{2+\frac{1}{2},1+\frac{1}{2}}
(K)$ for any compact subset $K$ of $M\times [0,T_4]$ to some function $\oa{F}$ 
as $k\to\infty$. Then $\oa{F}\in C^{2+\frac{1}{2},1+\frac{1}{2}}
(M\times [0,T_4])$ is a solution of (2.2) in $M\times (0,T_4)$ which 
satisfies (2.11). By Corollary 1.5 and an argument similar to the proof of 
Theorem 7.1 of Chapter VII of \cite{LSU} $\oa{F}\in
C^{3,1+\frac{1}{2}}(\2{B}_k\times (0,T_3])$ and the theorem follows.
\enddemo

\proclaim{\bf Theorem 2.6}
Let 
$$
T_1=\min (\log 2/(2C_3K_1^2),T)\tag 2.12
$$ 
where $C_3$ is given by Lemma 1.6.
Then for any $a\ge\max (1,3\pi/4\sqrt{k_0})$ there exists a solution 
$\oa{F}\in C^{2+\frac{1}{2},1+\frac{1}{2}}(M\times [0,T_1])\cap
C^{3,1+\frac{1}{2}}(M\times (0,T_1])$ of (2.2) 
in $M\times (0,T_1)$ which satisfies 
$$
|\oa{\nabla}\oa{F}|\le 2K_1\quad\text{ on }M\times (0,T_1).\tag 2.13
$$ 
\endproclaim
\demo{Proof}
Let $a\ge\max (1,3\pi/4\sqrt{k_0})$. Let $T_4$ and $C_6$ be as in 
Lemma 2.4. By Theorem 2.5 there exists a solution $\oa{F}\in 
C^{2+\frac{1}{2},1+\frac{1}{2}}(M\times [0,T_4])\cap C^{3,1+\frac{1}{2}}
(M\times (0,T_4])$ of (2.2) in $M\times (0,T_4)$ which satisfies (2.11). 
For any $R_1\ge 1$, let
$$
u(x,t)=|\oa{\nabla}\oa{F}(x,t)|^2\4{\phi}(x)
$$
where $\4{\phi}(x)=\phi(\rho(x)/R_1)$ with $\rho (x)=\rho (p_0,x)$. 
By the results in Chapter 1 of \cite{SY} $\rho (p_0,\cdot)\in C^{\infty}
(M\setminus (\text{Cut}(p_0)\cup\{p_0\}))$ (cf. \cite{GW1},\cite{W1},
\cite{W2}). Hence $u(\cdot,t)\in C^{\infty}(M\setminus\text{Cut}(p_0))$ 
for any $0<t\le T_4$. We first suppose that $x$ is not a cut point of 
$p_0$. Then by (2.9),
$$\align
u_t-\oa{\Delta}_tu
=&\4{\phi}\biggl (\frac{\1}{\1 t}-\oa{\Delta}_t\biggr)
|\oa{\nabla}\oa{F}|^2
-2\oa{\nabla}(|\oa{\nabla}\oa{F}|^2)\cdot\oa{\nabla}\4{\phi}
-|\oa{\nabla}\oa{F}|^2\oa{\Delta}_t\4{\phi}\\
\le&(-2|\oa{\nabla}^2\oa{F}|^2+C_3|\oa{\nabla}\oa{F}|^4)\4{\phi}
+(2/R_1)|\phi'||\oa{\nabla}(|\oa{\nabla}\oa{F}|^2)||\oa{\nabla}\rho|_{\oa{g}}\\
&\qquad +|\oa{\nabla}\oa{F}|^2\biggl(|\phi''|
\frac{|\oa{\nabla}\rho|_{\oa{g}}^2}{R_1^2}
+|\phi'|\frac{\oa{\Delta}_t\rho}{R_1}\biggr).\tag 2.14
\endalign
$$
Let $c_2>0$ be as given by Lemma 1.2. Then by (1.7),
$$
|\oa{\nabla}\rho|_{\oa{g}}^2=e^{-\psi_a}g^{ij}(x,t)\frac{\1\rho}{\1 x^i}
\frac{\1\rho}{\1 x^j}\le c_2g^{ij}(x,0)\frac{\1\rho}{\1 x^i}
\frac{\1\rho}{\1 x^j}\le c_2.\tag 2.15
$$
By (0.13) and (2.15),
$$
\frac{2}{R_1}|\phi'||\oa{\nabla}(|\oa{\nabla}\oa{F}|^2)|
|\oa{\nabla}\rho|_{\oa{g}}
\le\frac{4\sqrt{c_2}}{R_1}|\phi'||\oa{\nabla}F||\oa{\nabla}^2F|
\le|\oa{\nabla}^2F|^2\4{\phi}+\frac{C_8}{R_1^2}|\oa{\nabla}F|^2
\tag 2.16
$$
where
$$
C_8=4c_2\sup_{\Bbb{R}}(\phi'{}^2/\phi).
$$
By (2.4),
$$\align
\oa{\Delta}_t\rho=&\oa{g}^{ij}\biggl\{\frac{\1^2\rho}{\1 x^i\1 x^j}
-\oa{\Gamma}_{ij}^k\nabla_k\rho\biggr\}\\
=&e^{-\psi_a}\Delta_t\rho-\frac{1}{2}e^{-\psi_a}g^{ij}\biggl\{\delta_i^k
\nabla_j\psi_a+\delta_j^k\nabla_i\psi_a
-g^{kl}g_{ij}\nabla_l\psi_a\biggr\}\nabla_k\rho\\
=&e^{-\psi_a}\biggl (\Delta_t\rho+\frac{n-2}{2}g^{ij}\nabla_i
\rho\nabla_j\psi_a\biggr )\tag 2.17
\endalign
$$
By an argument similar to the proof of Lemma 1.3,
$$
\Delta_t\rho\le nc_2\frac{1+\sqrt{k_0}\rho}{\rho}+C_9\sqrt{t}\tag 2.18
$$
for some constant $C_9>0$ independent of $a$.
By Lemma 1.2 and Corollary 1.5,
$$
|g^{ij}\nabla_i\rho\nabla_j\psi_a|\le|\nabla\rho|_{g(t)}|\nabla\psi_a|_{g(t)}
\le c_2|\nabla\rho|_{g(0)}|\4{\nabla}\psi_a|_h\le C_{10}\tag 2.19
$$
for some constant $C_{10}>0$. By (2.14), (2.15), (2.16), (2.17), (2.18)
and (2.19),
$$
u_t-\oa{\Delta}_tu
\le (-|\oa{\nabla}^2\oa{F}|^2+C_3|\oa{\nabla}\oa{F}|^4)\4{\phi}
+\frac{C_{11}}{R_1}|\oa{\nabla}\oa{F}|^2
\le C_3|\oa{\nabla}\oa{F}|^2u+\frac{C_{11}}{R_1}|\oa{\nabla}\oa{F}|^2
\tag 2.20
$$
in $M\times (0,T_4)$ for some constant $C_{11}>0$. By (2.11) and (2.20),
$$\align
&u_t-\oa{\Delta}_tu\le C_3C_6^2u+\frac{C_{11}}{R_1^2}C_6^2
\qquad\qquad\qquad\qquad\qquad\quad\text{ in }
(M\setminus\text{Cut}(p_0))\times (0,T_4]\\
\Rightarrow\quad&\biggl(\frac{\1}{\1 t}-\oa{\Delta}_t\biggr)(e^{-C_3C_6^2t}u)
\le\frac{C_{11}C_6^2}{R_1}e^{-C_3C_6^2t}<\frac{2C_{11}C_6^2}{R_1}
\quad\text{ in }(M\setminus\text{Cut}(p_0))\times (0,T_4].\tag 2.21
\endalign
$$
Let
$$
q(y,t)=e^{-C_3C_6^2t}u(x,t)-(2C_{11}C_6^2/R_1)t.
$$
Then by (2.21),
$$
q_t<\oa{\Delta}_tq\tag 2.22
$$
in $(M\setminus\text{Cut}(p_0))\times (0,T_4]$. Suppose there exists
$(x_0,t_0)\in B(R_1)\times (0,T_4]$ such that
$$
q(x_0,t_0)=\max_{\2{B(R_1)}\times [0,T_4]}q(x,t).
$$
Suppose first that $x_0\in M\setminus\text{Cut}(p_0)$. Then
$$\align
&\frac{\1^2q}{\1x^i\1x^j}(x_0,t_0)\le 0, \frac{\1q}{\1 t}(x_0,t_0)\ge 0,
\frac{\1q}{\1 x^i}(x_0,t_0)=0\quad\forall i,j=1,2,\cdots,n\\
\Rightarrow\quad&q_t(x_0,t_0)\ge\oa{\Delta}_tq(x_0,t_0)\tag 2.23
\endalign
$$
which contradict (2.22). Hence $x_0\in \text{Cut}(p_0)$. Let $\gamma$ be a 
minimal geodesic in $(M,g(0))$ joining $p_0$ and $x_0$. Let 
$$
0<\3<\min (\rho (p_0,x_0),R_1/5).
$$
We choose a point $x_{\3}$ along the geodesic $\gamma$ such that 
$\rho(p_0,x_{\3})=\3$. Then $x_0$ is not a cut point 
of $x_{\3}$. Hence there exist a constant $\delta>0$ such that 
$B(x_0,\delta)\subset B(p_0,R_1)\setminus\text{Cut}(x_{\3})$.
Let 
$$
u_{\3}(x,t)=|\oa{\nabla}\oa{F}(x,t)|^2\phi((\3+\rho(x_{\3},x))/R_1)
\quad\forall x\in M, 0\le t\le T_4.
$$
and
$$
q_{\3}(x,t)=e^{-C_3C_6^2t}u_{\3}(x,t)-(2C_{11}C_6^2/R_1)t.
$$
By (2.21) and a similar argument as before $q_{\3}$ satisfies (2.22) in 
$B(x_0,\delta)\times (0,T_4]$.
Since
$$
\rho(p_0,x)\le\3+\rho(x_{\3},x)\quad\text{ in }B(x_0,\delta)
$$
and
$$
\rho(p_0,x_0)=\rho(p_0,x_{\3})+\rho(x_{\3},x_0)=\3+\rho(x_{\3},x_0),
$$
we have
$$
q_{\3}(x,t)\le q(x,t)\quad\text{ in }B(x_0,\delta)\times (0,T_4]
$$
and
$$
q_{\3}(x_0,t_0)=q(x_0,t_0).
$$
Hence $q_{\3}$ attains its maximum in $B_h(x_0,\delta)\times (0,T_4]$ at 
$(x_0,t_0)$. Thus $q_{\3}$ also satisfies (2.23). This contradicts (2.22)
for $q_{\3}$. Hence no such interior maximum point $(x_0,t_0)$ exists and
$$\align
&q(x,t)\le\max_{\1_pQ_{R_1}^{T_4}}q\le\sup_M|\nabla f|^2=K_1^2\qquad
\qquad\qquad\qquad\qquad\text{ in }Q_{R_1}^{T_4}\quad\forall R_1\ge 1\\
\Rightarrow\quad&e^{-C_3C_6^2t}|\oa{\nabla}\oa{F}(x,t)|^2
\phi(\rho(p_0,x)/R_1)\le K_1^2+(2C_{11}C_6^2/R_1)t\quad\text{ in }
Q_{R_1}^{T_4}\quad\forall R_1\ge 1\\
\Rightarrow\quad&|\oa{\nabla}\oa{F}(x,t)|\le K_1e^{\frac{C_3C_6^2}{2}t}
\qquad\qquad\quad\forall x\in M, 0\le t\le T_4\quad\text{ as }R_1\to\infty.
\tag 2.24
\endalign
$$
Let $T'=\min (2\log 2/C_3C_6^2,T_4)$. By (2.24), (2.13) holds in $M\times 
[0,T']$. Let 
$$\align
T_1'
=&\sup\{T''\in (0,T):\exists \text{ a solution of (2.2) in $M\times 
[0,T'')$ such that (2.13) holds}\\
&\qquad \text{ in } M\times (0,T'')\}.\endalign
$$ 
Then $T_1'\ge T'$. We claim that $T_1'\ge T_1$. Suppose not. Then 
$T_1'<T_1$. Since (2.13) holds in $M\times [0,T_1']$, by an argument 
similar to the proof of (2.24) but with (2.13) replacing (2.11) in the 
proof,
$$
|\oa{\nabla}\oa{F}(x,t)|\le K_1e^{\frac{C_3(2K_1)^2}{2}t}
=K_1e^{2C_3K_1^2t}\le K_1e^{2C_3K_1^2T_1'}<2K_1
\quad\forall x\in M, 0\le t\le T_1'.\tag 2.25
$$
By Theorem 2.5 there exists a solution $\oa{F}_1\in C^{2+\frac{1}{2},
1+\frac{1}{2}}(M\times [0,T_0])\cap C^{3,1+\frac{1}{2}}(M\times (0,T_0])$ 
of (2.2) in $M\times (0,T_0)$ with initial value $\oa{F}(x,T_1')$
for some constant $T_0\in (0,T-T_1')$ which satisfies
$$
|\oa{\nabla}\oa{F}_1|\le\4{C}_6\quad\text{ on }M\times [0,T_0]
$$
for some constant $\4{C}_6>0$. We extend $\oa{F}$ to a solution of (2.1) in 
$M\times (0,T_1+T_0)$ by setting $\oa{F}(x,t)=F_1(x,t-T_1')$ for any $x\in M$,
$T_1'\le t\le T_1'+T_0$. By (2.25) and an argument similar to the proof 
of (2.24),
$$
|\oa{\nabla}\oa{F}_1(x,t)|\le e^{\frac{C_3\4{C}_6^2}{2}t}
|\oa{\nabla}\oa{F}(x,T_1')|\le K_1e^{2C_3K_1^2T_1'}e^{\frac{C_3\4{C}_6^2}{2}t}
\quad\forall x\in M, 0\le t\le T_0\tag 2.26
$$
Let $\delta_2=-2C_3K_1^2T_1'+\log 2$. Then $\delta_2>0$.
Let $T_0'=\min (2\delta_2/(C_3\4{C}_6^2),T_0)$. Then by (2.26),
$$\align
&|\oa{\nabla}\oa{F}_1(x,t)|\le 2K_1\quad\text{ on }M\times [0,T_0]\\
\Rightarrow\quad&|\oa{\nabla}\oa{F}(x,t)|\le 2K_1\quad\text{ on }M\times 
[0,T_1'+T_0].
\endalign
$$
This contradicts the maximality of $T_1'$. Hence $T_1'\ge T_1$ and the 
theorem follows.
\enddemo

\proclaim{Theorem 2.7}
Let $f:(M,g(0))\to (N,h)$ be a given diffeomorphism satisfying (0.10) and 
(0.11) and let $T_1$ be given by (2.12). Then there exists a smooth solution
$F(\cdot,t):(M,g(t))\to (N,h)$ to the following Ricci harmonic flow:
$$\left\{\aligned
&\frac{\1F}{\1 t}=\Delta_{g(t),h}F\quad\text{ in }M\times (0, T_1)\\
&F(x,0)=f(x)\quad\text{ in }M
\endaligned\right.\tag 2.27
$$
which satisfies 
$$
|\nabla F|\le 2K_1\quad\text{ in }M\times [0,T_1]\tag 2.28 
$$
and
$$
\sup_{M\times [0,T_1)}|\nabla^mF|\le C_mt^{-\frac{m-2}{2}}\quad
\quad\text{ in }M\times [0,T_1]\quad\forall m\ge 2\tag 2.29
$$
for some constants $C_m>0$ depending on $k_0$, $K_1$ and $K_2$.
\endproclaim
\demo{Proof}
Let $a_1=\max (1,3\pi/4\sqrt{k_0})$. For any $a\ge a_1$, let 
$\oa{F}\in C^{2+\frac{1}{2},1+\frac{1}{2}}(M\times [0,T_1])\cap
C^{3,1+\frac{1}{2}}(M\times (0,T_1])$ be the solution of (2.2) in 
$M\times (0,T_1)$ given by Theorem 2.6 which satisfies (2.13).
By an argument similar to the proof of Theorem 7.1 of Chapter VII of 
\cite{LSU} (cf. P.245--246 of \cite{S1}) $\{\oa{F}\}_{a\ge a_1}$ is 
uniformly bounded in $C^{2+\frac{1}{2},1+\frac{1}{2}}(K)$ for any
compact subset $K$ of $M\times [0,T_1]$.  
By the Ascoli Theorem and a diagonalization argument the sequence 
$\{\oa{F}\}_{a\ge a_1}$ has a convergent subsequence 
$\{\ov{a_i}{F}\}_{i=1}^{\infty}$ that converges uniformly in 
$C^{2+\frac{1}{2},1+\frac{1}{2}}(K)$ for any compact subset $K$ of 
$M\times [0,T_1]$ to some function $F$ 
as $a_i\to\infty$. Then $F\in C^{2+\frac{1}{2},1+\frac{1}{2}}
(M\times [0,T_1])$.

Note that since $a\ge\max (1,3\pi/4\sqrt{k_0})$, by (1.13) and (1.25), 
$$
\2{\rho}_a(y)=0\quad\forall\rho(p_0,y)\le\frac{a}{2}-k_1
\quad\Rightarrow\quad\psi_a(y)=1\quad\forall\rho(p_0,y)\le\frac{a}{2}-k_1.
$$
Thus $F$ is a solution of (2.27) in $M\times (0,T_1)$. By (2.13) $F$
satisfies (2.28). 
By a bootrap argument and an argument similar to the proof of Theorem 7.1 
of Chapter VII of \cite{LSU} and that of P.245--246 of \cite{S1}
$F\in C^{\infty}(M\times (0,T_1))$. By (2.28) and an argument similar
to the proof of Lemma 2.12 of \cite{CZ} but with $\phi^a\equiv 0$ and
(2.28) replacing Lemma 2.11 in the proof there we get (2.29) and the 
theorem follows.
\enddemo

Note that by (0.8) and Section 6 of \cite{H3},
$$
\frac{\1}{\1 t}|\nabla F|^2\le\Delta_t|\nabla F|^2
-2|\nabla^2F|^2+k_0|\nabla F|^4\tag 2.30
$$
holds for any solution of (2.27). Hence by (2.30), Theorem 2.7 and an 
argument similar to the proof of Theorem 2.6 we have the following theorem.

\proclaim{Theorem 2.8}
Let $f:(M,g(0))\to (N,h)$ be a given diffeomorphism satisfying (0.10) and 
(0.11) and let $T_1=\min (\log 2/(2k_0K_1^2),T)$.
Then (2.27) has a smooth solution $F(\cdot,t):(M,g(t))\to (N,h)$ in
$M\times (0,T_1)$ which satisfies (2.28) and (2.29).
\endproclaim

Let $I:(M, g(0))\to (N,h)$ be the identity map. 
By direct computation,
$$
|\nabla I|_{g(0),h}=(g^{ij}(x,0)g_{\alpha\beta}(x,0)
\delta_i^{\alpha}\delta_j^{\beta})^{\frac{1}{2}}=\sqrt{n}
$$
and
$$
|\nabla^2I|_{g(0),h}=0.
$$
Hence by Theorem 2.8 we have the following result. 

\proclaim{Theorem 2.9}
Let 
$$
T_1=\min (\log 2/(2k_0n),T).
$$
Then there exists a smooth solution $F(\cdot,t):(M,g(t))\to (N,h)$ to the 
following Ricci harmonic flow:
$$\left\{\aligned
&\frac{\1F}{\1 t}=\Delta_{g(t),h}F\quad\text{ in }M\times (0, T_1)\\
&F(x,0)=x\qquad\quad\text{in }M
\endaligned\right.\tag 2.31
$$
which satisfies (2.28) and (2.29) with $K_1=\sqrt{n}$ for some constants 
$C_m>0$ depending on $k_0$.
\endproclaim

$$
\text{Section 3}
$$

In this section we will prove the uniqueness of solutions of Ricci flow 
on complete noncompact manifolds with bounded curvature. 

\proclaim{\bf Lemma 3.1}(cf. Proposition 3.1 of \cite{CZ})
Let $F$ and $T_1$ be as given by Theorem 2.9. Let $h=g(0)$ and let 
$\2{h}(t)=F^{\ast}(h)=(F(\cdot,t))^{\ast}(h)$ be the pull back metric 
of $h$ on $M$ by $F$. Then there exist a constant $C_1>0$ and a constant 
$0<T_2\le T_1$ depending only on $k_0$ such that
$$
\frac{1}{C_1}\2{h}_{ij}\le g_{ij}\le C_1\2{h}_{ij}
\quad\text{ in }M\times [0,T_2]\tag 3.1
$$
and
$$
|\2{\nabla}^kg|_{\2{h}}\le C_kt^{-\frac{k-1}{2}}\quad\text{ in }M\times 
[0,T_2]\quad\forall k\in\Bbb{Z}^+
\tag 3.2
$$
for some constants $C_k$ where $\2{\nabla}$ is the Levi-Civita connection
on $M$ with respect to $\2{h}$.
\endproclaim
\demo{Proof}
(3.2) and the first inequality on the the left hand side of (3.1) is
proved in Proposition 3.1 of \cite{CZ}. However the proof of the second 
inequality on the right hand side of (3.1) in \cite{CZ} is questionable. 
For the sake of completeness we will give a correct proof here.

By direct computation,
$$\align
\frac{\1\2{h}_{ij}}{\1 t}
=&\nabla_t(h_{\alpha\beta}\nabla_iF^{\alpha}\nabla_jF^{\beta})
=h_{\alpha\beta}(\nabla_t\nabla_iF^{\alpha})(\nabla_jF^{\beta})
+h_{\alpha\beta}(\nabla_iF^{\alpha})(\nabla_t\nabla_jF^{\beta})\\
=&h_{\alpha\beta}\nabla_i(F^{\alpha})_t\nabla_jF^{\beta}
+h_{\alpha\beta}\nabla_iF^{\alpha}\nabla_j(F^{\beta})_t.\tag 3.3
\endalign
$$
Now
$$\align
&h_{\alpha\beta}\nabla_i(F^{\alpha})_t\nabla_jF^{\beta}\\
=&h_{\alpha\beta}\nabla_i(\Delta_{g(t),h}F^{\alpha})\nabla_jF^{\beta}\\
=&h_{\alpha\beta}g^{i'j'}(\nabla_i\nabla_{i'}\nabla_{j'}F^{\alpha})
(\nabla_jF^{\beta})\\
=&h_{\alpha\beta}g^{i'j'}(\nabla_{i'}\nabla_i\nabla_{j'}F^{\alpha}
-R_{i'ij'l}g^{lk}\nabla_kF^{\alpha}+\4{R}_{pqrs}h^{\alpha r}
\nabla_iF^p\nabla_{i'}F^q\nabla_{j'}F^s)(\nabla_jF^{\beta})\\
=&h_{\alpha\beta}\Delta_t(\nabla_iF^{\alpha})\nabla_jF^{\beta}
-h_{\alpha\beta}g^{lk}R_{ik}\nabla_lF^{\alpha}\nabla_jF^{\beta}
+\4{R}_{\alpha\beta\gamma\delta}g^{kl}\nabla_jF^{\alpha}\nabla_kF^{\beta}
\nabla_iF^{\gamma}\nabla_lF^{\delta}.\tag 3.4
\endalign
$$
By (3.3) and (3.4),
$$\align
\frac{\1\2{h}_{ij}}{\1 t}
=&h_{\alpha\beta}\nabla_iF^{\alpha}\Delta_t
(\nabla_jF^{\beta})+h_{\alpha\beta}\Delta_t(\nabla_iF^{\alpha})
\nabla_jF^{\beta}
-h_{\alpha\beta}g^{lk}R_{ik}\nabla_lF^{\alpha}\nabla_jF^{\beta}\\
&\qquad -h_{\alpha\beta}g^{lk}R_{jk}\nabla_lF^{\alpha}\nabla_iF^{\beta}
+2g^{lk}\4{R}_{\alpha\beta\gamma\delta}\nabla_iF^{\alpha}
\nabla_kF^{\beta}\nabla_jF^{\gamma}\nabla_lF^{\delta}.\tag 3.5
\endalign
$$
Hence by (3.5), (2.29), (2.30) and Theorem 2.9, 
$$\align
&\biggl|\frac{\1\2{h}_{ij}}{\1 t}\biggr |_{g(t)}
\le C\sup_{M}(|\nabla F||\nabla^3F|+|\nabla F|^2+|\nabla F|^4)
\le\frac{C}{\sqrt{t}}\quad\text{ in }M\times (0,T_1)\\
\Rightarrow\quad&-\frac{C_1}{\sqrt{t}}g_{ij}\le\frac{\1\2{h}_{ij}}{\1 t}
\le\frac{C_1}{\sqrt{t}}g_{ij}\quad\text{ in }M\times (0,T_1)\tag 3.6
\endalign
$$
for some constant $C_1>0$. By Lemma 1.2 there exists constants $c_1,c_2>0$
such that (1.7) holds. By (1.7) and (3.6),
$$\align
\frac{\1\2{h}_{ij}}{\1 t}(x,t)\ge&-\frac{c_2C_1}{\sqrt{t}}g_{ij}(x,0)
\quad\text{ in }M\times (0,T)\\
\Rightarrow\quad\2{h}_{ij}(x,t)\ge&h_{ij}(x,0)-2c_2C_1\sqrt{t}g_{ij}(x,0)
\quad\text{ in }M\times (0,T)\\
=&g_{ij}(x,0)-2c_2C_1\sqrt{t}g_{ij}(x,0)\quad\text{ in }M\times (0,T)\\
\ge&(1-2c_2C_1\sqrt{t})g_{ij}(x,0)\quad\text{ in }M\times (0,T).\tag 3.7
\endalign
$$
Let $T_2=\min (T_1,1/(16c_2^2C_1^2))$. Then by (3.7) and Lemma 1.2,
$$
\2{h}_{ij}(x,t)\ge g_{ij}(x,0)/2\ge (1/2c_2)g_{ij}(x,t)\quad
\forall x\in M, 0\le t\le T_2
$$
and the lemma follows.
\enddemo

By Theorem 2.9, Lemma 3.1, and an argument similar to the proof of 
Proposition 3.2 of \cite{CZ} we have the following result.

\proclaim{\bf Lemma 3.2}(cf. Proposition 3.1 and Proposition 3.2 of \cite{CZ})
Let $F$ and $T_1$ be as given by Theorem 2.9. Then there exists a constant 
$0<T_2=T_2(g)\le T_1$ depending only on $k_0$ such that 
$F(\cdot,t):(M,g(t))\to (N,h)$ is a diffeomorphism for any $0\le t<T_2$. 
Moreover if $\hat{g}(t)$ is given by (0.4), then there exists a constant 
$C_1>0$ such that
$$
\frac{1}{C_1}h_{\alpha\beta}\le\hat{g}_{\alpha\beta}\le  C_1h_{\alpha\beta}
\quad\text{ in }N\times [0,T_2]\tag 3.8
$$
and
$$
|\hat{g}|_h^2+|\4{\nabla}\hat{g}|_h^2+t|\4{\nabla}^2\hat{g}|_h^2
\le C_1\quad\text{ in }N\times [0,T_2].
\tag 3.9
$$
\endproclaim

\proclaim{\bf Lemma 3.3}(cf. Proposition 3.3 of \cite{CZ})
Suppose $g(t)$ and $\2{g}(t)$ are both solutions of (0.1) in $M\times (0,T)$
with $g(0)\equiv \2{g}(0)$ on $M$ and both $Rm (g(t))$ and $Rm (\2{g}(t))$
satisfy (0.8) for some constant $k_0>0$. Let $F$ and $\2{F}$ be the solutions 
of (2.31) in $M\times (0,T_1(g))$ and in $M\times (0,T_1(\2{g}))$ 
for some constants $0<T_1(g),T_1(\2{g})\le T$  given by Theorem 2.9 
corresponding to Ricci flows $g(t)$ and $\2{g}(t)$ respectively.
Let $T_2(g)$ and $T_2(\2{g})$ be the constants given by Lemma 3.2.  
Let $\hat{g}(t)=(F(\cdot,t)^{-1})^{\ast}(g(t))$ and $\hat{\2{g}}(t)
=(F(\cdot,t)^{-1})^{\ast}(\2{g}(t))$ be the push forward metric of 
$g(t)$ and $\2{g}(t)$ on $N$ by $F$. Then $\hat{g}(t)\equiv\hat{\2{g}}(t)$ on
$M\times (0,T_0)$ where $T_0=\min (T_2(g),T_2(\2{g}))$.
\endproclaim
\demo{Proof}
A proof of this lemma (Proposition 3.3 of \cite{CZ}) is given in \cite{CZ}. 
However the proof given in \cite{CZ} is not correct because the 
deduction of the last two inequalities on P.151 of \cite{CZ} assumed 
that one can interchange differentiation and taking limit as $\3\to 0$
which is not true in general.
For the sake of completeness we will modify their argument and give a correct 
proof of the result here. We will use the technique of proof of Theorem 2.6
to proof this lemma. We first recall that by the proof on P.150-151 of 
\cite{CZ}, we have
$$
\frac{\1}{\1 t}|\hat{g}-\hat{\2{g}}|_h^2\le\hat{g}^{\alpha\beta}
\4{\nabla}_{\alpha}\4{\nabla}_{\beta}|\hat{g}-\hat{\2{g}}|_h^2
+\frac{C_2}{\sqrt{t}}|\hat{g}-\hat{\2{g}}|_h^2\quad\text{ in }N\times (0,T_0)
\tag 3.10
$$
for some constant $C_2>0$ where $h=g(0)=\2{g}(0)$
and $\4{\nabla}$ is the covariant derivative of $h$. For any $R_1\ge 1$ we let
$$
u(y,t)=|\hat{g}-\hat{\2{g}}|_h^2\4{\phi}(y)\quad\forall y\in N, 
0\le t\le T_0.
$$
where $\4{\phi}(y)=\phi(\rho(y)/R_1)$ and $\rho(y)=\rho(p_0,y)$. 
We first suppose that $y$ is not a cut point of $p_0$. Then by (3.8), 
(3.9), (3.10), and Lemma 1.2 for any $0<t\le T_0$,
$$\align
&u_t-\hat{g}^{\alpha\beta}\4{\nabla}_{\alpha}\4{\nabla}_{\beta}u\\
=&\biggl [\biggl (\frac{\1}{\1 t}-\hat{g}^{\alpha\beta}
\4{\nabla}_{\alpha}\4{\nabla}_{\beta}\biggr )
|\hat{g}-\hat{\2{g}}|_h^2\biggr ]\4{\phi}
-2\hat{g}^{\alpha\beta}\4{\nabla}_{\alpha}|\hat{g}-\hat{\2{g}}|_h^2
\cdot\4{\nabla}_{\beta}\4{\phi}-(\hat{g}^{\alpha\beta}
\4{\nabla}_{\alpha}\4{\nabla}_{\beta}\4{\phi})
|\hat{g}-\hat{\2{g}}|_h^2\\
\le&\frac{C_2}{\sqrt{t}}|\hat{g}-\hat{\2{g}}|_h^2\4{\phi}
+(4/R_1)|\hat{g}|_h|\hat{g}-\hat{\2{g}}|_h|\4{\nabla}_{\alpha}\hat{g}
-\4{\nabla}_{\alpha}\hat{\2{g}}|_h|\phi'||\4{\nabla}\rho|\\
&\qquad+|\hat{g}-\hat{\2{g}}|_h^2
\biggl(\frac{|\phi'|\hat{g}^{\alpha\beta}\4{\nabla}_{\alpha}\4{\nabla}_{\beta}
\rho}{R_1}+\frac{|\phi''||\hat{g}^{\alpha\beta}||\4{\nabla}\rho|^2}{R_1^2}
\biggr )\\
\le&C_2\frac{u}{\sqrt{t}}+\frac{C_3'}{R_1}+\frac{2C_3'}{R_1}|\phi'|
\hat{g}^{\alpha\beta}\4{\nabla}_{\alpha}\4{\nabla}_{\beta}\rho
\tag 3.11
\endalign
$$
for some constant $C_3'>0$. By the Hessian comparison theorem \cite{SY},
P.309-310 of \cite{S2}, and (3.8),
$$\align
&\4{\nabla}_{\alpha}\4{\nabla}_{\beta}\rho
\le\frac{1+\sqrt{k_0}\rho}{\rho}h_{\alpha\beta}\\
\Rightarrow\quad&\hat{g}^{\alpha\beta}\4{\nabla}_{\alpha}
\4{\nabla}_{\beta}\rho\le\hat{g}^{\alpha\beta}h_{\alpha\beta}
\frac{1+\sqrt{k_0}\rho}{\rho}\le nC_1\frac{1+\sqrt{k_0}\rho}{\rho}.
\tag 3.12
\endalign
$$
By (3.11) and (3.12),
$$\align
&u_t-\hat{g}^{\alpha\beta}\4{\nabla}_{\alpha}\4{\nabla}_{\beta}u
\le C_2\frac{u}{\sqrt{t}}+\frac{C_4}{R_1}\qquad\qquad\qquad\qquad
\qquad\text{in }(N\setminus\text{Cut}_h(y_0))\times (0,T_0]\\
\Rightarrow\quad&\biggl (\frac{\1}{\1 t}-\hat{g}^{\alpha\beta}
\4{\nabla}_{\alpha}\4{\nabla}_{\beta}\biggr )(ue^{-2C_2\sqrt{t}})\le
\frac{C_4}{R_1}e^{-2C_2\sqrt{t}}<\frac{2C_4}{R_1}\quad\text{ in }
(N\setminus\text{Cut}_h(y_0))\times (0,T_0].\tag 3.13
\endalign
$$
for some constant $C_4>0$. Let
$$
q(y,t)=u(y,t)e^{-2C_2\sqrt{t}}-(2C_4/R_1)t.
$$
Then by (3.13),
$$
q_t<\hat{g}^{\alpha\beta}\4{\nabla}_{\alpha}\4{\nabla}_{\beta}q
\quad\text{ in }(N\setminus\text{Cut}_h(p_0))\times (0,T_0].\tag 3.14
$$
By (3.14) and an argument
similar to the proof of Theorem 2.6 the function $q$ attains its maximum 
on $\1_p(B_h(y_1,\delta)\times (0,T_0))$. Hence 
$$\align
&q(y,t)\le\max_{\1_p(B_h(y_0,R_1)\times [0,T_0])}q(y,t)=0
\qquad\qquad\,\,\,\text{in }B_h(y_0,R_1)\times [0,T_0]\\
\Rightarrow\quad&|\hat{g}-\hat{\2{g}}|_h^2\phi((\rho(p_0,y))/R_1)
\le(2C_4/R_1)te^{2C_2\sqrt{t}}\quad\text{ in }
\2{B_h(p_0,R_1)}\times [0,T_0]\quad\forall R_1>1\\
\Rightarrow\quad&|\hat{g}-\hat{\2{g}}|_h^2=0\qquad\qquad\qquad\qquad\qquad
\qquad\qquad\quad\text{ in }N\times [0,T_0]\text{ as }R_1\to\infty
\endalign
$$
and the lemma follows.
\enddemo

By the same argument as the proof on P.152 of \cite{CZ} but with Lemma 3.4
replacing Proposition 3.3 in the proof there we get the following uniqueness
theorem.

\proclaim{\bf Theorem 3.5}
Suppose $g(t)$ and $\2{g}(t)$ are both solutions of (0.1) in $M\times (0,T)$
with $g(0)\equiv \2{g}(0)$ on $M$ and both $Rm (g(t))$ and $Rm (\2{g}(t))$
satisfy (0.8) for some constant $k_0>0$. Then $g(t)\equiv\2{g}(t)$ on
$M\times (0,T)$.
\endproclaim

\enddocument